\documentclass[12pt]{article}
\usepackage[numbers]{natbib}
\usepackage{array,epsfig,fancyheadings,rotating}

\usepackage{amsmath}
\usepackage{amssymb}
\usepackage{amsfonts}
\usepackage{amsthm}

\usepackage{graphicx}

\usepackage{bbm}
\usepackage{float}
\usepackage{sectsty}

\sectionfont{\large}
\subsectionfont{\normalsize}

\RequirePackage[OT1]{fontenc}
\RequirePackage{amsthm,amsmath}
\RequirePackage{natbib}

\numberwithin{equation}{section}
\newtheorem{thm}{Theorem}
\newtheorem{lm}{Lemma}

\newtheorem{pro}{Proposition}
\theoremstyle{definition}

\newtheorem{re}{Remark}

\begin{document}


\centerline{\large\bf INTERMEDIATE EFFICIENCY OF TESTS} 
\centerline{\large\bf UNDER HEAVY-TAILED ALTERNATIVES}
\vspace{.4cm} \centerline{Tadeusz Inglot} \vspace{.4cm} \centerline{\it
Wroc{\l}aw University of Science and Technology} 


\begin{quotation}
\noindent {\bf Abstract.}
  We show that for local alternatives which are not square integrable the intermediate (or Kallenberg) efficiency of the Neyman-Pearson test for uniformity with respect to the classical Kolmogorov-Smirnov test is equal to infinity. Contrary to this, for local square integrable alternatives the intermediate efficiency is finite and can be explicitly calculated. \par

\vspace{9pt}
\par

\vspace{2mm}
\end{quotation}\par


\section{Introduction and testing problem}\label{S1}

We consider the classical problem of testing for uniformity. We compare the Neyman-Pearson (NP) test with the classical Kolmogorov-Smirnov (KS) test for uniformity for a class of local unbounded alternatives in terms of asymptotic relatve efficiency (ARE) notion. By ARE we mean the Kallenberg's intermediate efficiency which is a limit of the ratio of sample sizes which guarantee the same precision for both tests (the same significance level tending to 0 slower than exponentially and the same asymptotically nondegenerate power). 

Our main issue is that for alternatives which are not square integrable the efficiency of the KS test with respect to the NP test cannot be positive. In particular, we apply the simplest variant of the intermediate efficiency notion recently elaborated in Inglot et al. \cite{r8} and called pathwise intermediate efficiency. We show that this efficiency for the NP test with respect to the KS test for a class of alternatives approaching the null distribution, which are not square integrable, is equal to $\infty$ (Theorem \ref{thm1}).

Recall that the notion of the intermediate efficiency was introduced originally by Kallenberg \cite{r9}. Then it was developed and applied to some testing problems and several tests in a series of papers in the last two decades e.g. Inglot \cite{r3}, Inglot and Ledwina \cite{r5,r6,r7}, Mason and Eubank \cite{r10}, Mirakhmedov \cite{r11} or recently Inglot et al. \cite{r8} and \'Cmiel et al. \cite{r1}. For more detailed discussions and up-to-date remarks and comments we send the reader to Inglot et al. \cite{r8}. Note that, by the definition, this efficiency notion involves asymmetric requirements for compared tests.

In Inglot and Ledwina \cite{r7} it was found, among others, the intermediate efficiency of the KS test with respect to the NP test  for sequences of bounded alternatives approaching the null distribution. In the present paper, as a byproduct, we extend that result to unbounded square integrable alternatives in the reversed formulation i.e. taking the KS test as a benchmark procedure and comparing the NP test to it.

Since we consider both simple testing problem as well as very regular statistics and to make the paper self-contained we do not refer to general results and technical tools elaborated in Inglot et al. \cite{r8}. Instead, we present all auxiliary results and all proofs directly.\\

Let $X_1,...,X_n$ be independent random variables with values in $[0,1]$ and a distribution $P$ with continuous distribution function. By $P_0$ we denote the uniform distribution over the interval $[0,1]$. We test the simple null hypothesis
$$ H_0: P=P_0$$ 
against 
$$ H_1: P\neq P_0.$$
To compare tests consider local alternatives with densities (with respect to $P_0$) of the form
\begin{equation}\label{1.1}
p_{\theta_n}(t)=1-\theta_n+\theta_n f(t),\;t\in(0,1),
\end{equation}
where $\theta_n\in (0,1),\;\theta_n\to 0$ as $n\to\infty$ and $f$ is a fixed alternative density.

By $P_{\theta_n}$ we denote the distribution with the density $p_{\theta_n}(t)$. Moreover, $P_0^n, \;P_{\theta_n}^n$ shall denote $n$ fold products of $P_0$ and $P_{\theta_n}$, respectively. 

For each $n$ consider the standardized NP test statistic
\begin{equation}\label{1.2} 
V_n=\frac{1}{\sqrt{n}\sigma_{0n}}\sum_{i=1}^n (\log p_{\theta_n}(X_i)-e_{0n})
\end{equation}
for testing $H_0$ against the simple hypothesis $H_{1n}: P=P_{\theta_n}$. Here
$$e_{0n}=\int_0^1\log p_{\theta_n}(t) dt,\;\;\;\sigma^2_{0n}=\int_0^1\log^2p_{\theta_n}(t)dt-e_{0n}^2$$
are the two first moments of $\log p_{\theta_n}(X_1)$ under $P_0$ which are finite due to the integrability of $f$. Additionally denote 
$$e_n=\int_0^1p_{\theta_n}(t)\log p_{\theta_n}(t)dt,\;\;\;\sigma_n^2=\int_0^1p_{\theta_n}(t)\log^2p_{\theta_n}(t)dt-e_n^2$$
the corresponding moments under $P_{\theta_n}$ which we assume to be finite. For example, this is the case if $f\in L_q[0,1]$ for some $q>1$. 

Now, set 
\begin{equation}\label{1.3} 
b_n=\sqrt{n}(e_n-e_{0n})/\sigma_{0n}.
\end{equation}
The sequence $b_n$ shall play a role of an asymptotic shift of $V_n$ under $P_{\theta_n}$. 

For each $n$ and any fixed $x\in \mathbb{R}$ set
\begin{equation}\label{1.4} 
\alpha_n=\alpha_n(x)=P_0^n(V_n\geq x+b_n)
\end{equation}
the significance level of the NP test corresponding to the critical value $x+b_n$. Since $V_n$ is bounded in probability under $P_0$, then whenever $b_n\to\infty$ we have $\alpha_n\to 0$.

Let
$$K_n=\sqrt{n}\sup_{t\in (0,1)}|\hat{F}_n(t)-t|,$$
where $\hat{F}_n(t)$ is the empirical distribution function of $X_1,...,X_n$, be the classical unweighted KS test statistic. For each $n$ and every $N\geq n$ let $u_{N,n}$ be the exact critical value of the KS test at the level $\alpha_n$ defined by (\ref{1.4}) and for the sample size $N$ i.e.
$$ P_0^N(K_N\geq u_{N,n})=\alpha_n.$$
For each $n$ let $N_n$ be the minimal sample size such that for all $k\geq 0$
\begin{equation}\label{1.5} 
P_{\theta_n}^{N_n+k}(K_{N_n+k}\geq u_{N_n+k,n})\geq P_{\theta_n}^n(V_n\geq x+b_n)
\end{equation}
i.e. the minimal sample size begining from which the power of the KS test under $P_{\theta_n}$ and at the level $\alpha_n$ is not smaller than that for the NP test at the same lavel and for the sample size $n$. Obviously, $N_n\geq n$. The limit of the ratio $N_n/n$, if exists, is called the intermediate efficiency of the NP test with respect to the KS test (cf. Inglot et al. \cite{r8}). We study an asymptotic behaviour of the ratio $N_n/n$ succesively for two cases when $f$ is heavy tailed or square integrable and show that they lead to qualitatively different answers.

The paper is organized as follows. In Section \ref{S2} we consider local alternatives which are not square integrable while in Section \ref{S3} square integrable ones. In Section \ref{S4} we present outcomes of a simulation study nicely illustrating theoretical results. All proofs are sent to Sections \ref{S5} -- \ref{S8}.

In the sequel we shall use the following notation: for sequences $x_n,y_n$ of positive numbers by $x_n\asymp y_n$ we shall mean that for some positive constants $c_1, c_2$ it holds $c_1\leq x_n/y_n\leq c_2$ for all $n$ while by $x_n\sim y_n$ we shall mean that $x_n/y_n\to 1$ as $n\to\infty$.

\section{Heavy tailed case}\label{S2}

Assume that a density $f$ in (\ref{1.1}) satisfies the following condition:\\

\hspace*{0.6cm}for some $r\in (0,1)$ and some positive $C_1\leq (1-r)^r$ and $C_2>1$ we have
\begin{equation}\label{2.1}
\quad\quad\quad C_1t^{-r}\leq f(t)\leq C_2t^{-r},\;\;\mbox{for}\;t\in (0,C_1^r),\;\;\mbox{and}\;\; f(t)\leq C_2\;\;\mbox{for}\; t\in [C_1^{1/r},1).
\end{equation}
Observe that $f\in L_2[0,1]$ if and only if $r\in (0,1/2)$. In the present section we consider the case $r\in[1/2,1)$.

First we describe an asymptotic behaviour of $b_n$, defined in (\ref{1.3}). It is an immediate corollary of Lemma \ref{le3} proved in Section \ref{S6}.\\

\begin{pro}\label{pro1} 
If $f$ satisfies (\ref{2.1}) for some $r\in[1/2,1)$ then
\begin{equation}\label{2.2}
b_n\asymp \sqrt{n}\kappa_{nr},
\end{equation}
where
$$\kappa_{nr}=\left\{\begin{array}{ll} \theta_n^{1/2r},&if\;\;r\in(1/2,1)\\\theta_n\sqrt{\log(1/\theta_n)},&if\;\;r=1/2\end{array}\right. .$$
\end{pro}

\vspace{2mm}
The next proposition is a simple consequence of (\ref{2.2}).

\begin{pro}\label{pro2} 
Let $p_{\theta_n}(t)$ be a sequence of densities given by (\ref{1.1}) with $f$ satisfying (\ref{2.1}) for some $r\in[1/2,1)$ and $\theta_n\to 0$ is such that $n\kappa_{nr}^2\to\infty$. Then for every $x\in \mathbb{R}$ it holds
\begin{equation}\label{2.3}
0<\liminf_{n\to\infty} P_{\theta_n}^n(V_n\geq x+b_n)\leq\limsup_{n\to\infty} P_{\theta_n}^n(V_n\geq x+b_n)<1.
\end{equation}
\end{pro}

\vspace{3mm}
\noindent The proof of Proposition \ref{pro2} is given in Section \ref{S6}.

\begin{thm}\label{thm1} 
Let $p_{\theta_n}(t)$ be a sequence of densities given by (\ref{1.1}) with $f$ satisfying (\ref{2.1}) for some $r\in[1/2,1)$ and $\theta_n\to 0$ is such that $n\kappa_{nr}^2\to\infty$. Then for any $x\in \mathbb{R}$ and the significance levels defined by (\ref{1.4}) we have for $N_n$ defined by (\ref{1.5})
\begin{equation}\label{2.4}
\lim_{n\to\infty}\frac{N_n}{n}=\infty.
\end{equation}
\end{thm}

The proof of Theorem \ref{thm1} is given in Section \ref{S5}.

\begin{re}\label{re1}
In terms of the intermediate efficiency (as defined in Inglot et al. \cite{r8}) Theorem \ref{thm1} says that for $f$ satisfying (\ref{2.1}) this efficiency of the NP test with respect to the KS test is equal to $\infty$. This efficiency notion requires a nondegenerate asymptotic power of the second compared test, here the NP test. It is essential in the proof of Theorem \ref{thm1} and is ensured by our Proposition \ref{pro2}. In the proof of Theorem \ref{thm1} we directly show that the intermediate slope of the KS test equals $2n\theta_n^2||A||_{\infty}^2$ without introducing such terminology and without referring to regularity conditions (I.1) and (I.2) in Inglot et al. \cite{r8}. For the NP test the regularity condition (II.2) (ibid.) can be deduced from the proofs of Proposition \ref{pro2} and Lemma \ref{le4}. Moreover, it is enough to show a weaker property than the regularity condition (II.1) (ibid.) meaning that an expression which may be considered as the intermediate slope of the NP test is at least of order $n\kappa_{nr}^2$. Anyway, we prove (\ref{2.4}) in the simplest possible way.
Obviously, the statement (\ref{2.4}) remains true for any test for uniformity which has positive and finite intermediate efficiency with respect to the KS test and simultaneously can be taken as a benchmark procedure. For some further comments see Section 2 in Inglot et al. \cite{r8}.
\end{re}

\begin{re}\label{re2}
The assumption that $f$ is unbounded at the left end of $(0,1)$ is not essential. Obviously, our result is valid for $f$ unbounded at the right end of $(0,1)$ or at both ends (not necesserily symmetrically) or in some interior point of $(0,1)$, as well, provided a condition analogous to (\ref{2.1}) is satisfied.
\end{re}

\begin{re}\label{re3}
In \'Cmiel et al. \cite{r1} the intermediate efficiency of some weighted goodness of fit tests has been investigated. In particular, from results of that paper it follows that, opposite the statement of Theorem \ref{thm1}, for $f\in L_q[0,1],\;q>1,$ the intermediate efficiency of the integral Anderson-Darling test with respect to the KS test is finite with an explicit formula for calculating it. Also, for $f\in L_q[0,1],\;q>2,$ the intermediate efficiencies of the classical Anderson-Darling (weighted supremum) test and its truncated version exist with an explicit formulae for them (cf. Remark 4, ibid.). 
\end{re}
\section{Square integrable case}\label{S3}

Suppose that $f$ in (\ref{1.1}) belongs to $L_2[0,1]$. Set
$$ a(t)=\frac{1}{c} (f(t)-1),$$
where $c^2=\int_0^1(f(t)-1)^2dt$. Then by rescaling $\theta_n$ we may rewrite (\ref{1.1}) in the equivalent form 
$$p_{\theta_n}(t)=1+\theta_n a(t),\;\;t\in (0,1).$$ 

In the present setting an asymptotic behaviour of $b_n$ in (\ref{1.3}), stated below, is an immediate corollary of Lemma \ref{le5} proved in Section \ref{S8}.

\begin{pro}\label{pro3} 
If $f\in L_2[0.1]$ then
\begin{equation}\label{3.1}
b_n\sim \sqrt{n}\theta_n\;\;\mbox{and}\;\;\sigma_{0n}\sim \theta_n.
\end{equation}
\end{pro}

A consequence of Proposition \ref{pro3} is the following result which plays the same role as Proposition \ref{pro2} in the heavy-tailed case.

\begin{pro}\label{pro4} 
Let $p_{\theta_n}(t)$ be a sequence of densities given by (1.1) with $f\in L_2[0,1]$ and $\theta_n\to 0$ is such that $n\theta_n^2\to\infty$. Then for every $x\in \mathbb{R}$ it holds
\begin{equation}\label{3.2}
0<\liminf_{n\to\infty} P_{\theta_n}^n(V_n\geq x+b_n)\leq\limsup_{n\to\infty} P_{\theta_n}^n(V_n\geq x+b_n)<1.
\end{equation}
\end{pro}

\noindent The proof of Proposition \ref{pro4} is given in Section \ref{S8}. Now, we state our second main result.

\begin{thm}\label{thm2} 
Let $p_{\theta_n}(t)$ be a sequence of densities given by (\ref{1.1}) with $f\in L_2[0,1]$ and $\theta_n\to 0$ is such that $n\theta_n^2\to\infty$. Then for any $x\in \mathbb{R}$ and the significance levels defined by (\ref{1.4}) we have for $N_n$ defind by (\ref{1.5})
\begin{equation}\label{3.3}
\lim_{n\to\infty}\frac{N_n}{n}=\frac{1}{4||A||_{\infty}^2}={\cal E}(a),
\end{equation}
where 
$$A(t)=\int_0^t a(u)du$$
and $||\cdot||_{\infty}$ denotes the supremum norm on $[0,1]$.
\end{thm}

\begin{re}\label{re4} 
Theorem \ref{thm2} says that the intermediate efficiency (as defined in Inglot et al. \cite{r8}), of the NP test with respect to the KS test for converging square integrable sequences of alternatives exists and equals $1/4||A||_{\infty}^2$. Thus it extends Corollary 6.2 of Inglot and Ledwina \cite{r7} to the case of unbounded square integrable alternatives. Note that Corollary 6.2 was stated equivalently in terms of the intermediate efficiency of the KS test with respect to the NP test. Note also that in the proof of Theorem \ref{thm2} we find the intermediate slopes of compared tests equal to $2n\theta_n^2||A||_{\infty}^2$ and $n\theta_n^2/2$, respectively, under the assumptions of this theorem without introducing such terminology.
\end{re}

The proof of Theorem \ref{thm2} is given in Section \ref{S7}.\\

\noindent{\bf Example.} For $r\in (0,1/2)$ let $f_r(t)=(1-r)t^{-r},\;t\in (0,1)$, and consequently $a_r(t)=(\sqrt{1-2r}/r)((1-r)t^{-r}-1)$. Then 
\begin{equation}\label{3.4}
{\cal E}(a_r)=\frac{(1-r)^{2-2/r}}{4(1-2r)}.
\end{equation}
Observe that ${\cal E}(a_r)\to \infty$ when $r\to 1/2$ which nicely agrees with the statement of Theorem \ref{thm1}. 


\section{Simulation results}\label{S4}

Below, we present results of a small simulation study showing how (\ref{2.4}) and (\ref{3.3}) are reflected empirically for a particular density 
$$f_r(t)=(1-r)t^{-r},\;t\in(0,1).$$ 
We select some small values of $\theta_n=\theta$ and keep powers separated from 0 and 1. We take heavy-tailed alternatives by choosing two values of $r$ greater than 1/2 and square integrable alternatives represented by two values of $r$ smaller than 1/2. In the two last cases the formula (3.4) can be appplied. The results are shown in Tables 1 -- 4.

\begin{center}
{\small {\sc Table 1.} Empirical powers (in \%) of the NP and KS tests for the\\ alternative $f_r$, small values of $\theta$ and several $n$. $\alpha=0.05,\;r=0.7$.}\\

\vspace{2mm}
\begin{tabular}{c|rr||c|rr||c||rr}
$n$&\multicolumn{2}{c||}{$\theta=0.1$}&$n$&\multicolumn{2}{c||}{$\theta=0.05$}&$n$&\multicolumn{2}{c}{$\theta=0.02$}\\ 
&KS&NP&&KS&NP&&KS&NP\\ \hline
16&5&30&11&4&15&40&4&15\\
28&6&40&20&4&20&70&5&20\\
41&7&50&42&5&30&155&5&30\\
60&9&60&70&6&40&250&5&40\\
80&11&70&105&7&50&3200&15&99\\
300&30&98&150&7&60&4900&20&100 \\
410&40&100&540&15&94&7700&30&100\\
520&50&100&750&20&98&10100&40&100\\
640&60&100&1200&30&100&&&\\
780&70&100&1600&40&100&&&\\
&&&2080&50&100&&&\\
&&&2500&60&100&&&\\ \hline
\end{tabular}
\end{center}

\vspace{5mm}
\begin{center}
{\small {\sc Table 2.} Empirical powers (in \%) of the NP and KS tests for the\\ alternative $f_r$, small values of $\theta$ and several $n$. $\alpha=0.05,\;r=0.6$.}\\

\vspace{2mm}
\begin{tabular}{c|rr||c|rr||c|rr}
$n$&\multicolumn{2}{c||}{$\theta=0.1$}&$n$&\multicolumn{2}{c||}{$\theta=0.05$}&$n$&\multicolumn{2}{c}{$\theta=0.02$}\\ 
&KS&NP&&KS&NP&&KS&NP\\ \hline
16&5&20&27&5&15&110&4&15\\
35&6&30&48&5&20&205&5&20\\
60&7&40&105&6&30&460&6&30\\
92&9&50&180&7&40&5400&15&95\\
127&10&60&260&8&50&7800&20&99\\
175&13&70&370&9&60&12000&30&100\\
300&20&86&800&15&85&&&\\
480&30&96&1045&20&94&&&\\
640&40&99&1900&30&99&&&\\
840&50&100&2600&40&100&&&\\
1040&60&100&3300&50&100&&&\\
1280&70&100&4100&60&100&&&\\ \hline
\end{tabular}\end{center}

\vspace{5mm}
\begin{center}
{\small {\sc Table 3.} Empirical powers (in \%) of the NP and KS tests for the alter-\\native $f_r$, small values of $\theta$ and several $n$. $\alpha=0.05,\;r=0.4, {\cal E}(a_{0.4})=5.787$.}\\

\vspace{2mm}
\begin{tabular}{c|rr||c|rr||c||rr}
$n$&\multicolumn{2}{c||}{$\theta=0.2$}&$n$&\multicolumn{2}{c||}{$\theta=0.1$}&$n$&\multicolumn{2}{c}{$\theta=0.05$}\\ 
&KS&NP&&KS&NP&&KS&NP\\ \hline
 15& 5&15&54&5&15&205&6&15\\
 28& 6&20&105&6&20&400&6&20\\
 62& 9&30&220&8&30&810&8&30\\
100&11&40&360&11&40&2400&15&59\\
148&14&50&510&13&50&3400&20&72\\
153&15&51&600&15&55&5600&30&88\\
200&18&60&700&16&60&&\\
225&20&64&870&20&68&&&\\
270&24&70&1430&30&84&&&\\
350&30&79&1950&40&93&&&\\
500&40&90&2500&50&97&&&\\
640&50&95&3160&60&99&&&\\ 
795&60&97&&&&&&\\
970&70&99&&&&&&\\ \hline
\end{tabular}
\end{center}

\vspace{1mm}
\begin{center}
{\small {\sc Table 4.} Empirical powers (in \%) of the NP and KS tests for the alter-\\native $f_r$, small values of $\theta$ and several $n$. $\alpha=0.05,\;r=0.3,\; {\cal E}(a_{0.3})=3.302$.}\\

\vspace{1mm}
\begin{tabular}{c|rr||c|rr||c||rr}
$n$&\multicolumn{2}{c||}{$\theta=0.2$}&$n$&\multicolumn{2}{c||}{$\theta=0.1$}&$n$&\multicolumn{2}{c}{$\theta=0.05$}\\ 
&KS&NP&&KS&NP&&KS&NP\\ \hline
40&6&15  &160&5&15  & 640&6&15\\
75&7&20  &300&7&20  &1200&7&20\\
165&10&30&610&10&30 &2300&10&30\\
255&13&40&950&13&40 &4800&15&47\\
300&15&44&1200&15&47&6800&20&59\\
360&17&50&1340&17&50&11100&30&77\\
430&20&56&1700&20&58&&&\\
480&22&60&1830&21&60&&&\\
645&28&70&2800&30&76&&&\\
710&30&73&3880&40&87&&&\\
980&40&84&5050&50&93&&&\\
1260&50&91&6350&60&97&&&\\ 
1600&60&96&&&&&&\\
1950&70&98&&&&&&\\ \hline
\end{tabular}
\end{center}

\vspace{4mm}
Using the results from Tables 1 -- 4 we present in Table 5 ratios $N_n/n$ for four considered values of $r$, some small values of $\theta$ and several powers separated from 0 and 1.

From Table 5 it is easily seen that for $r>1/2$ the ratio $N_n/n$ behaves unstably and rapidly grows when $\theta$ tends to 0 thus confirming the statement of Theorem \ref{thm1}. Contrary to this, for $\;r<1/2\;$ the ratio behaves stably and takes values relatively close to the 

\noindent intermediate efficiency of the NP test with respect to the KS test given by the formula (\ref{3.4}).

\vspace{2mm}
\begin{center}
{\small {\sc Table 5.} Ratios of $N_n/n$ for the alternative $f_r$, small values of $\theta$, several\\ powers separated from 0 and 1 and four values of $r$. Significance level $\alpha=0.05$.}\\

\vspace{2mm}
\begin{tabular}{c|c|ccccccc}
&&\multicolumn{7}{c}{\small power in \%}\\
$r$&$\theta$&15&20&30&40&50&60&70\\ \hline
0.7&0.10&&&          18.8&14.6&12.7&10.7&9.8\\
   &0.05&49.1&37.5&28.6&22.9&19.8&16.7\\
	 &0.02&80.0&70.0&49.7&40.4&&&\\ \hline
0.6&0.10&&18.8&13.7&10.7&9.1&8.2&7.3\\
&   0.05&29.6&21.8&18.1&14.4&12.7&11.1&\\
&   0.02&49.1&38.1&26.1&&&&\\ \hline
0.4&0.20&10.2    &8.0&5.7&5.0&4.3&4.0&3.6\\
   &0.10&11.1&8.3&6.6&5.4&4.9&4.5&\\
	 &0.05&11.7&8.5&6.9&&&&\\ \hline
0.3&0.20&7.5&5.7&4.3&3.8&3.5&3.3&3.0\\
   &0.10&7.5&5.7&4.6&4.1&3.8&3.5&\\
	 &0.05&7.5&5.7&4.8&&&&\\ \hline
\end{tabular}
\end{center}


\section{Proof of Theorem \ref{thm1}}\label{S5}

A key step in the proof of our theorem is a moderate deviation result both for $V_n$ and $K_n$ under the null distribution. Below we state it as two separate propositions. The first one is stated in a weak version but sufficient to prove Theorem \ref{thm1}.

\begin{pro}\label{pro5} 
If $f$ in (\ref{1.1}) satisfies (\ref{2.1}) for some $r\in [1/2,1)$ then for every sequence $x_n$ of positive numbers such that $x_n=O(\kappa_{nr})$ we have
\begin{equation}\label{5.1}
-\limsup_{n\to\infty}\frac{1}{nx_n^2}\log P_0^n(V_n\geq \sqrt{n}x_n)>0.
\end{equation}
\end{pro}

The proof of Proposition \ref{pro5} is given in Section \ref{S6}. The pertaining moderate deviation theorem for $K_n$ was obtained in Inglot and Ledwina \cite{r4}. For completeness we state it below.

\begin{pro}\label{pro6}
For every sequence $x_n$ of positive numbers such that $x_n\to 0$ and $n x_n^2\to\infty$ it holds
\begin{equation}\label{5.2}
-\lim_{n\to\infty}\frac{1}{nx_n^2}\log P_0^n(K_n\geq \sqrt{n}x_n)=2.
\end{equation}
\end{pro}

\vspace{3mm}
Now, we are ready to prove the theorem. Take any $x\in \mathbb{R}$. Proposition \ref{pro2} says that the sequence of powers of the NP test at the significance level $\alpha_n$ defined by (\ref{1.4}) is bounded away from 0 and 1. Set $x_n=(x+b_n)/\sqrt{n}$. Then by Proposition \ref{pro1} $x_n\asymp \kappa_{nr}$ and from Proposition \ref{pro5} it follows that for some positive constants $c, c'$ and sufficiently large $n$
\begin{equation}\label{5.3}
-\log\alpha_n=-\log P_0^n(V_n\geq x+b_n)\geq c'nx_n^2\geq cn\kappa_{nr}^2.
\end{equation}
Set $A(t)=\int_0^tf(u)du-t$ and $F_n(t)=t+\theta_nA(t)$. Then by the triangle inequality and for $N_n$ defined by (\ref{1.5}) we have
$$P_{\theta_n}^{N_n}(K_{N_n}\geq u_{N_n,n})=Pr(||e_{N_n}\circ F_n+\sqrt{N_n}\theta_nA||_{\infty}\geq u_{N_n,n})$$
$$\leq Pr(||e_{N_n}||_{\infty}\geq u_{N_n,n}-\sqrt{N_n}\theta_n||A||_{\infty}),$$
where $e_N(t)$ denotes the uniform empirical process for the sample of size $N$ while $Pr$ denotes a probability on the underlying probability space. From (\ref{1.5}), Proposition \ref{pro2} and the convergence of $e_{N_n}$ in distribution to a Brownian bridge it follows that for some positive $C$ 
$$u_{N_n,n}-\sqrt{N_n}\theta_n||A||_{\infty}\leq C.$$
This implies $u_{N_n,n}/\sqrt{N_n}\to 0$. Since $P_0^{N_n}(K_{N_n}\geq u_{N_n,n})=\alpha_n$ and $\alpha_n\to 0$, then $u_{N_n,n}\to\infty$ and Proposition \ref{pro6} applied to $x_n=u_{N_n,n}/\sqrt{N_n}$ gives
$$ -\log \alpha_n=2u^2_{N_n,n}(1+o(1)).$$
This together with (\ref{5.3}) gives for sufficiently large $n$
$$cn\kappa_{nr}^2\leq 2u_{N_n,n}^2(1+o(1))\leq 2(C+\sqrt{N_n}\theta_n||A||_{\infty})^2(1+o(1))\leq 5C^2+5N_n\theta_n^2||A||_{\infty}^2.$$
As $\theta_n/\kappa_{nr}\to 0$ the above implies $n/N_n\to 0$ and finishes the proof of (\ref{2.4}). \hfill $\Box$
\section{Proofs of Propositions \ref{pro1}, \ref{pro2} and \ref{pro5}}\label{S6}

\subsection{Auxiliary lemmas}\label{S6.1}

For $k=0,1$ and integer $m\geq 1$ consider the following integrals 
$$ I_{km}(n)=\int_0^1 [\theta_ng(t)]^{k}\log^m(1+\theta_ng(t))dt,$$
$$ J_{km}(n)=\int_0^1 [1+\theta_ng(t)]^{k}|\log(1+\theta_ng(t))-e_{0n}|^mdt,$$
where $f$ and $\theta_n$ are as in (\ref{1.1}) and, for short, we have denoted $g(t)=f(t)-1$. The first lemma describes an asymptotic behaviour of $I_{km}(n)$ and $J_{km}(n)$ as $n\to\infty$ under $r>1/2$.

\begin{lm}\label{le1} 
Suppose $f$ satisfies (\ref{2.1}) for some $r\in[1/2,1)$. Then for any $k=0,1$ and any integer $m\geq 1$ such that $k+m\geq 2$ we have
\begin{equation}\label{6.1}
 I_{km}(n)\asymp \kappa_{nr}^2.
\end{equation}
Moreover, for any $k=0,1$ and $m\geq 2$ we have
\begin{equation}\label{6.2}
J_{km}(n)\asymp \kappa_{nr}^2.
\end{equation}
\end{lm}

The proof of Lemma \ref{le1} is based on the following elementary fact.

\begin{lm}\label{le2} 
Suppose $f$ satisfies (\ref{2.1}) for some $r\in (0,1)$. Then for any $k=0,1$ and any integer $m\geq 1$ we have
\begin{equation}\label{6.3}
\int_0^{C_1^{1/r}} [\theta_ng(t)]^{k}\log^m(1+\theta_ng(t))dt \asymp \left\{\begin{array}{ll}\theta_n^{\min\{k+m,1/r\}},&if\;\;k+m\neq 1/r\\\theta_n^{1/r}\log(1/\theta_n),&if\;\;k+m=1/r\end{array}\right.
\end{equation}
and
\begin{equation}\label{6.4}
\int_0^{C_1^{1/r}} [1+\theta_ng(t)]^{k}|\log(1+\theta_ng(t))-e_{0n}|^mdt \asymp \left\{\begin{array}{ll}\theta_n^{\min\{m,1/r\}},&if\;\;m\neq 1/r\\\theta_n^{1/r}\log(1/\theta_n),&if\;\;m=1/r\end{array}\right. .
\end{equation}
\end{lm}

\noindent{\bf Proof of Lemma 2.} When $k+m$ is odd then the function $\psi_{km}(y)=y^k\log^m(1+y)$ is increasing on $(-1,\infty)$ while for $k+m$ even $\psi_{km}(y)$ is decreasing on $(-1,0)$ and increasing on $(0,\infty)$. The condition (\ref{2.1}) implies that $f(t)\geq 1$ for $t\in(0,C_1^{1/r}]$. Hence, by the monotonicity of $\psi_{km}(y)$ on $(0,\infty)$ and the inequality $y/2\leq \log(1+y)\leq y$ holding on $(0,1/2)$, from (\ref{2.1}) and after the substitution $y=\theta_n(C_1/t^r-1)$, the integral in (\ref{6.3}) can be estimated for $n$ sufficiently large from below by 
$$\int_0^{C_1^{1/r}}[\theta_n(C_1/t^r-1)]^{k}\log^m(1+\theta_n(C_1/t^r-1))dt=r^{-1}C_1^{1/r}\theta_n^{1/r}\int_0^{\infty}y^{k}\frac{\log^m(1+y)}{(\theta_n+y)^{1+1/r}}dy$$
$$\geq 2^{-1-1/r}r^{-1}C_1^{1/r}\left[\theta_n^{1/r}\int_{1/2}^{\infty}y^{k}\frac{\log^m(1+y)}{y^{1+1/r}}dy + 2^{-m}\theta_n^{1/r}\int_{\theta_n}^{1/2}y^{k+m-1-1/r}dy\right]$$
and, after the substitution $y=\theta_n(C_2/t^r-1)$, for $n$ sufficiently large from above by
$$\int_0^{C_1^{1/r}}[\theta_n(C_2/t^r-1)]^{k}\log^m(1+\theta_n(C_2/t^r-1))dt=r^{-1}C_2^{1/r}\theta_n^{1/r}\int_{(C_2/C_1-1)\theta_n}^{\infty} y^k\frac{\log^m(1+y)}{(\theta_n+y)^{1+1/r}}dy$$
$$\leq r^{-1}C_2^{1/r}\left[\theta_n^{1/r}\int_{1/2}^{\infty} y^{k}\frac{\log^m(1+y)}{y^{1+1/r}}dy+\theta_n^{1/r}\int_{(C_2/C_1-1)\theta_n}^{1/2}y^{k+m-1-1/r}dy\right].$$
Since the second terms in the above estimates are of order $\theta_n^{\min\{k+m,1/r\}}$ if $k+m\neq 1/r$ or $\theta_n^{1/r}\log(1/\theta_n)$ if $k+m=1/r$, the relation (\ref{6.3}) is proved.

Now, observe that for $n$ sufficiently large
\begin{equation}\label{6.5}
0=\theta_n\int_0^1(f(t)-1)dt\leq -e_{0n}=-\int_0^1\log p_{\theta_n}(t) dt\leq-\log(1-\theta_n)\leq 2\theta_n.
\end{equation}
To prove (\ref{6.4}) we argue similarly as for (\ref{6.3}). Using (\ref{2.1}) and (\ref{6.5}) we estimate the integral in (\ref{6.4}) for sufficiently large $n$ from below by
$$\int_0^{C_1^{1/r}}[1+\theta_n(C_1/t^r-1)]^k(\log(1+\theta_n(C_1/t^r-1))-e_{0n})^mdt$$
$$= r^{-1}C_1^{1/r}\theta_n^{1/r}\int_0^{\infty}(1+y)^k\frac{(\log(1+y)-e_{0n})^m}{(\theta_n+y)^{1+1/r}}dy$$
$$\geq 2^{-1-1/r}r^{-1}C_1^{1/r}\left[\theta_n^{1/r}\int_{1/2}^{\infty}\frac{\log^m(1+y)}{y^{1+1/r}}dy+2^{-m}\theta_n^{1/r}\int_{\theta_n}^{1/2}y^{m-1-1/r}dy\right]$$
and from above by
$$\int_0^{C_1^{1/r}}[1+\theta_n(C_2/t^r-1)]^k(\log(1+\theta_n(C_2/t^r-1))-e_{0n})^mdt$$
$$=r^{-1}C_2^{1/r}\theta_n^{1/r}\int_{(C_2/C_1-1)\theta_n}^{\infty}(1+y)^k\frac{(\log(1+y)-e_{0n})^m}{(\theta_n+y)^{1+1/r}}dy$$
\begin{equation}\label{6.6}
\leq 2^mr^{-1}C_2^{1/r}\left[\theta_n^{1/r}\int_{1/2}^{\infty}(1+y)^k\frac{\log^m(1+y)}{y^{1+1/r}}dy+\left(\frac{3}{2}\right)^k\theta_n^{1/r}\int_{(C_2/C_1-1)\theta_n}^{1/2}(\theta_n+y)^{m-1-1/r}dy\right].
\end{equation}
Since the second terms in the above estimates are of order $\theta_n^{\min\{m,1/r\}}$ if $m\neq 1/r$ or $\theta_n^{1/r}\log(1/\theta_n)$ if $m=1/r$, the relation (\ref{6.4}) is proved. \hfill $\Box$\\

\noindent
{\bf Proof of Lemma 1.}  The monotonicity properties of the functions $\psi_{km}(y)$ defined in the proof of Lemma \ref{le2} and (\ref{2.1}) imply that for $n$ sufficiently large
$$ \left|I_{km}(n)-\int_0^{C_1^{1/r}}[\theta_n(f(t)-1)]^k\log^m(1+\theta_n(f(t)-1))dt\right|$$
$$\leq \theta_n^k|\log(1-\theta_n)|^m+[\theta_nC_2]^k\log^m(1+\theta_nC_2)\asymp \theta_n^{k+m}$$
and due to (\ref{6.5})
$$\left|J_{km}(n)-\int_0^{C_1^{1/r}}[1+\theta_n(f(t)-1)]^k(\log(1+\theta_n(f(t)-1))-e_{0n})^mdt\right|$$
$$\leq (1+\theta_nC_2)^k(\log(1+\theta_nC_2)-e_{0n})^m\leq (1+\theta_nC_2)^k[\theta_n(C_2+2)]^m\asymp \theta_n^{m}.$$
Since $r>1/2$ and $k+m\geq 2$, (\ref{6.1}) follows from (\ref{6.3}) while (\ref{6.2}) follows from (\ref{6.4}) and the assumption $m\geq 2$. This completes the proof of Lemma \ref{le1}.\hfill $\Box$

\begin{lm}\label{le3} If $f$ satisfies (\ref{2.1}) for some $r\in[1/2,1)$ then
\begin{equation}\label{6.7} 
e_n-e_{0n}\asymp \kappa_{nr}^2,\;\;\;\mbox{and}\;\;\; \sigma_{0n}^2\asymp \kappa_{nr}^2\asymp\sigma_n^2.
\end{equation}
\end{lm}

\vspace{3mm}\noindent
{\bf Proof.} Observe that $e_n-e_{0n}=I_{11}(n)$, $\sigma_{0n}^2=I_{02}(n)-e_{0n}^2$ and 
$\sigma_n^2=J_{12}(n)-(e_n-e_{0n})^2.$ 
Hence, (\ref{6.7}) follows immediately from Lemma \ref{le1}. \hfill $\Box$

\begin{lm}\label{le4} For each $n\geq 1$ let $X_1, X_2,...,X_n$ be independent random variables with density $p_{\theta_n}(t)$ and $f$ satisfying (\ref{2.1}) with $r\in(1/2,1)$. If $\theta_n\to 0$ such that $n\kappa_{nr}^2\to\infty$ then for every $y\in \mathbb{R}$
$$\lim_{n\to\infty} P_{\theta_n}^n\left(\frac{1}{\sqrt{n}\sigma_n}\sum_{i=1}^n(\log p_{\theta_n}(X_i)-e_n)\leq y\right)=\Phi(y),$$
where $\Phi(y)$ denotes the standard normal distribution function.
\end{lm}

\noindent
{\bf Proof.} Denote $Y_{ni}=\log p_{\theta_n}(X_i)-e_n,\;i=1,...,n,\;n\geq 1,$ the triangular array of independent mean 0 random variables. To prove Lemma \ref{le4} it is enough to check the Liapunov condition. We have
$E_{\theta_n}|Y_{ni}|^3\leq 4J_{13}(n)+4(e_n-e_{0n})^3\asymp \kappa_{nr}^2$ by (\ref{6.2}) and (\ref{6.7}). Since $\sigma_n^3\asymp \kappa_{nr}^3$ by (\ref{6.7}), the Liapunov condition holds true due to the assumption $n\kappa_{nr}^2\to\infty$.\hfill $\Box$

\subsection{Proof of Proposition \ref{pro2}}\label{S6.2}

Observe that
$$ V_n=\frac{\sigma_n}{\sigma_{0n}}\left[\frac{1}{\sqrt{n}\sigma_n}\sum_{i=1}^n(\log p_{\theta_n}(X_i)-e_n)\right]+b_n.$$
So, for $x\in \mathbb{R}$
$$P_{\theta_n}^n(V_n\geq x+b_n)=P_{\theta_n}^n\left(\frac{1}{\sqrt{n}\sigma_n}\sum_{i=1}^n(\log p_{\theta_n}(X_i)-e_n)\geq x\frac{\sigma_{0n}}{\sigma_n}\right)$$
and (\ref{2.3}) is an immediate consequence of Lemma \ref{le3} and Lemma \ref{le4}.\hfill $\Box$

\subsection{Proof of Proposition \ref{pro5}}\label{S6.3}

We shall apply the following version of the Bernstein inequality (cf. Yurinskii \cite{r12}).\\

\noindent
{\bf Theorem A.} {\it Let $\xi_1,...,\xi_n,\;n\geq 1,$ be independent identically distributed random variables with $E\xi_1=0$ and $E\xi_1^2=1$ such that for some constant $M>0$  it holds
\begin{equation}\label{6.8}
E|\xi_1|^m\leq \frac{m!}{2}M^{m-2}\;\;for\; every \;\;m\geq 3.
\end{equation}
Then for all $x>0$}
\begin{equation}\label{6.9}
P\left(\frac{\xi_1+...+\xi_n}{\sqrt{n}}\geq x\right)\leq 2\exp\left\{-\frac{x^2}{2(1+xM/\sqrt{n})}\right\}.
\end{equation}

\vspace{3mm}
In Theorem A set 
$$\xi_i=\frac{\log p_{\theta_n}(X_i)-e_{0n}}{\sigma_{0n}},\;\;i=1,...,n,$$ 
where $X_1,...,X_n$ are uniformly distributed over $[0,1]$. Then $E_0|\xi_1|^m=J_{0m}/\sigma_{0n}^m$ for all $m\geq 3$. Since
$$ \int_{1/2}^{\infty}\frac{\log^m(1+y)}{y^{1+1/r}}dy\leq 3^m\int_0^{\infty}\frac{\log^m(1+y)}{(1+y)^{1+1/r}}dy=r(3r)^m m!$$
then from (\ref{6.6}) it follows that (\ref{6.8}) holds with e.g. $M_n=6r/\sigma_{0n}$. Applying (\ref{6.9}) to $x=\sqrt{n}x_n$ we get
$$ P_0^n(V_n\geq \sqrt{n}x_n)\leq 2\exp\left\{-\frac{nx_n^2}{2(1+x_nM_n)}\right\}.$$
By the assumption and Lemma \ref{le3} we have $x_nM_n=O(\kappa_{nr}/\sigma_{0n})=O(1)$ and hence (\ref{5.1}) follows. \hfill $\Box$

\section{Proof of Theorem \ref{thm2}}\label{S7}

We shall apply the following moderate deviation result for $V_n$, proved in Section \ref{S8}.

\begin{pro}\label{pro7} If $f$ in (1.1) satisfies $f\in L_2[0,1]$ and $n\theta_n^2\to\infty$, then for any positive $\delta<1/2$ and every sequence $x_n$ satisfying $2\delta\sigma_{0n}<x_n<2(1-\delta)\sigma_{0n}$ it holds
\begin{equation}\label{7.1}
-\lim_{n\to\infty}\frac{1}{nx_n^2}\log P_0^n(V_n\geq \sqrt{n}x_n)=\frac{1}{2}.
\end{equation}
\end{pro}

\vspace{3mm}
Now, take any $x\in \mathbb{R}$. Proposition \ref{pro4} says that the sequence of powers of the NP test at the significance level $\alpha_n$ defined by (\ref{1.4}) is bounded away from 0 and 1. Set $x_n=(x+b_n)/\sqrt{n}$. Then, by Proposition \ref{pro3}, $x_n$ satisfies the assumption of Proposition \ref{pro7} for sufficiently large $n$. Hence
\begin{equation}\label{7.2}
-\log\alpha_n=-\log P_0^n(V_n\geq x+b_n) =\frac{n\theta_n^2}{2}(1+o(1)).
\end{equation}
Set $F_n(t)=t+\theta_nA(t)$. Recall that here $A(t)$, defined in Theorem \ref{thm2}, corresponds to the normalized function $a$. Then by the triangle inequality and $N_n$ defined in (\ref{1.5}) we have
$$P_{\theta_n}^{N_n}(K_{N_n}\geq u_{N_n,n})=Pr(||e_{N_n}\circ F_n+\sqrt{N_n}\theta_nA||_{\infty}\geq u_{N_n,n})$$
$$\leq Pr(||e_{N_n}||_{\infty}\geq u_{N_n,n}-\sqrt{N_n}\theta_n||A||_{\infty}),$$
where $e_N(t)$ denotes the uniform empirical process for the sample of size $N$. From (\ref{1.5}), Proposition \ref{pro4} and the convergence of $e_{N_n}$ in distribution to a Brownian bridge it follows that for some positive $C$ 
\begin{equation}\label{7.3}
u_{N_n,n}-\sqrt{N_n}\theta_n||A||_{\infty}\leq C.
\end{equation}
This implies $u_{N_n,n}/\sqrt{N_n}\to 0$. Since $P_0^{N_n}(K_{N_n}\geq u_{N_n,n})=\alpha_n$ and $\alpha_n\to 0$, then $u_{N_n,n}\to\infty$ and Proposition \ref{pro6} applied to $x_n=u_{N_n,n}/\sqrt{N_n}$ gives
$$ -\log \alpha_n=2u^2_{N_n,n}(1+o(1)).$$
This together with (\ref{7.2}) gives
$$\frac{n\theta_n^2}{2}= 2u_{N_n,n}^2(1+o(1))\leq 2(C+\sqrt{N_n}\theta_n||A||_{\infty})^2(1+o(1))$$
which means that
\begin{equation}\label{7.4}
\limsup_{n\to\infty}\frac{n}{N_n}\leq 4||A||_{\infty}^2.
\end{equation}
On the other hand, using the minimality property of $N_n$ in (\ref{1.5}) and Proposition \ref{pro4}, a similar argument as used to get (\ref{7.3}) leads to the relation
\begin{equation}\label{7.5}
u_{N_n-1,n}-\sqrt{N_n-1}\theta_n||A||_{\infty}\geq -C
\end{equation}
for some positive constant $C$. Observe that $u_{N_n-1,n}/\sqrt{N_n-1}\to 0$. Indeed, Proposition \ref{pro6} applies to $x'_n=\theta_n\sqrt{n/(N_n-1)}\to 0$ and gives $\log P_0^{N_n-1}(K_{N_n-1}\geq x'_n\sqrt{N_n-1})=-2n\theta_n^2(1+o(1))$ which together with (\ref{7.2}) and the definition of $u_{N_n-1,n}$ imply for $n$ sufficiently large $u_{N_n-1,n}/\sqrt{N_n-1}\leq x'_n$ thus proving our claim. Again applying Proposition \ref{pro6} to $x_n=u_{N_n-1,n}/\sqrt{N_n-1}$ we obtain 
$$ -\log \alpha_n=2u^2_{N_n-1,n}(1+o(1))$$
and consequently from (\ref{7.2}) and (\ref{7.5})
$$\frac{n\theta_n^2}{2}= 2u_{N_n-1,n}^2(1+o(1))\geq 2(\sqrt{N_n-1}\theta_n||A||_{\infty}-C)^2.$$
Hence,
$$\liminf_{n\to\infty}\frac{n}{N_n}\geq 4||A||_{\infty}^2$$
wich together with (\ref{7.4}) proves (\ref{3.3}). \hfill $\Box$

\section{Proofs of Propositions \ref{pro3}, \ref{pro4} and \ref{pro7}}\label{S8}

Recall some useful simple inequalities
\begin{equation}\label{8.1} 
\log^2(1+y)\leq y,\;\;y\geq 0,
\end{equation}
\begin{equation}\label{8.2}
\log^3(1+y)\leq \min\left\{\frac{3}{2}y,y^2\right\},\;\;y\geq 0,
\end{equation}
and for any $0<\varepsilon<1/2$
\begin{equation}\label{8.3}
(1-\varepsilon)y^2\leq y\log(1+y)\leq (1+\varepsilon)y^2,\;\;y\in[-\varepsilon,\varepsilon],
\end{equation}
\begin{equation}\label{8.4}
(1-\varepsilon)y^2\leq \log^2(1+y)\leq (1+2\varepsilon)y^2,\;\;y\in[-\varepsilon,\varepsilon].
\end{equation}
Proposition \ref{pro3} is an immediate corollary of the following lemma.

\begin{lm}\label{le5} If $f\in L_2[0,1]$ then $e_n-e_{0n}=\theta_n^2(1+o(1)),\;\; \sigma_{0n}^2=\theta_n^2(1+o(1))$ and $\sigma_n^2=\theta_n^2(1+o(1))$.
\end{lm}

\noindent
{\bf Proof.} Taking in (\ref{8.3}) $\varepsilon=\sqrt{\theta_n}$ and remembering that $a(t)\geq -1 \;a.s.$ we have for sufficiently large $n$
$$e_n-e_{0n}=\int_0^1\theta_na(t)\log(1+\theta_na(t))dt\geq \int_{a<1/\sqrt{\theta_n}}\theta_na(t)\log(1+\theta_na(t))dt$$
$$\geq (1-\sqrt{\theta_n})\theta_n^2\int_{a<1/\sqrt{\theta_n}}a^2(t)dt=\theta_n^2(1+o(1))$$
and
$$e_n-e_{0n}=\int_{a<1/\sqrt{\theta_n}}\theta_na(t)\log(1+\theta_na(t))dt+\int_{a\geq 1/\sqrt{\theta_n}}\theta_na(t)\log(1+\theta_na(t))dt$$
$$\leq (1+\sqrt{\theta_n})\theta_n^2\int_{a<1/\sqrt{\theta_n}}a^2(t)dt+\theta_n^2\int_{a\geq 1/\sqrt{\theta_n}}a^2(t)dt=\theta_n^2(1+o(1))$$ 
which proves the first statement.

From (\ref{8.4}) with $\varepsilon=\sqrt{\theta_n}$ and similar estimates as above we get
\begin{equation}\label{8.5}
\int_0^1\log^2(1+\theta_na(t))dt=\theta_n^2(1+o(1)).
\end{equation}
Moreover, from an obvious inequality $y-y^2\leq\log(1+y)\leq y,\;y\in[-1/2,1/2],$ it follows
\begin{equation}\label{8.6}
0\geq e_{0n}\geq \theta_n\int_{a\geq 1/\sqrt{\theta_n}}a(t)dt-\theta_n^2\int_{a< 1/\sqrt{\theta_n}}a^2(t)dt=o(\theta_n).\end{equation}
Combining (\ref{8.5}) and (\ref{8.6}) gets the second statement.

To prove the third one, note first that $e_n=(e_n-e_{0n})+e_{0n}=o(\theta_n)$. Moreover, since $a(t)\geq -1\;a.s.$ and by (\ref{8.1})
$$\left|\int_0^1\theta_na(t)\log^2(1+\theta_na(t))dt\right|$$
$$\leq\int_{a<0} \theta_n|a(t)|\log^2(1+\theta_na(t))dt+\int_{0<a<1/\sqrt{\theta_n}}\theta_na(t)\log^2(1+\theta_na(t))dt+\theta_n^2\int_{a>1/\sqrt{\theta_n}}a^2(t)dt$$
$$\leq\theta_n^3\int_{a<0}a^2(t)dt+\theta_n^{5/2}\int_{0<a<1/\sqrt{\theta_n}}a^2(t)dt+\theta_n^2\int_{a\geq 1/\sqrt{\theta_n}}a^2(t)dt=o(\theta_n^2).$$
Hence, from (\ref{8.5}) and the above
$$\sigma^2_n=\int_0^1\log^2(1+\theta_na(t))dt+\int_0^1\theta_na(t)\log^2(1+\theta_na(t))dt-e_n^2=\theta_n^2(1+o(1)).$$
This completes the proof of Lemma \ref{le5}. \hfill $\Box$

\begin{lm}\label{le6} For each $n\geq 1$ let $X_1, X_2,...,X_n$ be independent random variables with density $p_{\theta_n}(t)$ given by (\ref{1.1}) with $f\in L_2[0,1]$ and $\theta_n\to 0$ such that $n\theta_n^2\to\infty$ then for every $y\in \mathbb{R}$
$$\lim_{n\to\infty} P_{\theta_n}^n\left(\frac{1}{\sqrt{n}\sigma_n}\sum_{i=1}^n(\log (1+\theta_na(X_i))-e_n)\leq y\right)=\Phi(y).$$
\end{lm}

\vspace{2mm}\noindent
{\bf Proof.} Denote $Y_{ni}=\log(1+\theta_na(X_i))-e_n,\;i=1,...,n,\;n\geq 1,$ the triangular array of independent mean 0 random variables. To prove Lemma \ref{le6} it is enough to check the Liapunov condition. We have from (\ref{8.2}), (\ref{8.6}) and Lemma \ref{le5} for sufficiently large $n$
$$E_{\theta_n}|Y_{ni}|^3=\int_0^1(1+\theta_na(t))|\log(1+\theta_na(t))-e_n|^3dt$$
$$\leq 4e_n^3+4\int_{a<1/\sqrt{\theta_n}}(1+\theta_na(t))|\log(1+\theta_na(t))|^3dt+4\int_{a\geq 1/\sqrt{\theta_n}}(1+\theta_na(t))\log^3(1+\theta_na(t))dt$$
$$\leq 4e_n^3+4\sqrt{\theta_n}\int_{a<1/\sqrt{\theta_n}}(1+\theta_na(t))\log^2(1+\theta_na(t))dt+10\theta_n^2\int_{a\geq 1/\sqrt{\theta_n}}a^2(t)dt$$
$$\leq 4e_n^3+4\sqrt{\theta_n}(\sigma_n^2+e_n^2)+10\theta_n^2\int_{a\geq 1/\sqrt{\theta_n}}a^2(t)dt=o(\theta_n^2).$$
Hence the Liapunov condition holds due to Lemma \ref{le5} and the assumption $n\theta_n^2\to\infty$.\hfill $\Box$\\

It is easily seen that applying Lemma \ref{le6} the proof of Proposition \ref{pro4} goes exactly in the same way as that of Proposition \ref{pro2}.\\

The proof of Proposition \ref{pro7} is based on the following moderate deviation result of Ermakov \cite{r2}.\\

\noindent
{\bf Theorem B.} {\it Let $Y_{n1},...Y_{nn},\;n\geq 1,$ be the triangular array of independent identically distributed random variables with $EY_{n1}=0,$ Var$\,Y_{n1}=1$ and for some sequence $h_n>0,\;h_n\to 0,\; nh_n^2\to\infty$ it holds

(i) $\displaystyle E e^{h_nY_{n1}}<C;$

(ii) $\displaystyle E|Y_{n1}|^3\leq C\frac{\omega_n}{h_n}$ for some $\omega_n>0$ (possibly dependent on $h_n$).

Then for all $x$ such that $\delta h_n\leq x\leq (1-\delta)h_n$ for some $\delta>0$ it holds}
\begin{equation}\label{8.7}
\log P\left(\frac{Y_{n1}+...Y_{nn}}{\sqrt{n}}\geq \sqrt{n}x\right)=-\frac{nx^2}{2}+O(nh_n^2\omega_n).
\end{equation}

\vspace{2mm}
In Theorem B set $Y_{ni}=(\log(1+\theta_na(X_i))-e_{0n})/\sigma_{0n}$ and $h_n=2\sigma_{0n}$. Then $nh_n^2\to\infty$ by the assumption and Lemma \ref{le5}. Moreover, for sufficiently large $n$
$$E_0e^{h_nY_{n1}}=e^{-2e_{0n}}\int_0^1(1+\theta_na(t))^2dt<2$$
which proves the condition (i) in this theorem. From (8.2) and (8.6) we have similarly as previously
$$E_0|Y_{n1}|^3=\frac{1}{\sigma_{0n}^3}\int_0^1|\log(1+\theta_na(t))-e_{0n}|^3dt$$
$$\leq \frac{4|e_{0n}|^3}{\sigma_{0n}^3}+\frac{4\sqrt{\theta_n}}{\sigma_{0n}^3}\int_{a<1/\sqrt{\theta_n}}\log^2(1+\theta_na(t))dt+\frac{4}{\sigma_{0n}^3}\int_{a\geq 1/\sqrt{\theta_n}}\log^3(1+\theta_na(t))dt$$
$$\leq \frac{4}{\sqrt{\theta_n}}(1+o(1))+\frac{4}{\theta_n}\int_{a\geq 1/\sqrt{\theta_n}}a^2(t)dt.$$
So, the condition (ii) of Theorem B holds with $\omega_n=\max\{\sqrt{\theta_n},\int_{a\geq 1/\sqrt{\theta_n}}a^2(t)dt\}$tending to 0 and (\ref{7.1}) follows from (\ref{8.7}) by inserting $x_n$ in place of $x$.\hfill $\Box$


\renewcommand\bibname

\vspace{8mm}
\noindent
Tadeusz Inglot\\
Faculty of Pure and Applied Mathematics,\\ Wroc{\l}aw University of Science and Technology,\\
Wybrze\.ze Wyspia\'nskiego 27, 50-370 Wroc{\l}aw,\\ Poland.
\vskip 2pt
\noindent
E-mail: Tadeusz.Inglot@pwr.edu.pl

\end{document}